\documentclass[12pt]{amsart}
\usepackage{amsthm,vmargin}
\usepackage{amsmath,xypic}
\usepackage{hyperref}
\begin{document}
\newcommand{\Q}{{\mathbb Q}}
\newcommand{\C}{{\mathbb C}}
\newcommand{\R}{{\mathbb R}}
\newcommand{\Z}{{\mathbb Z}}
\newcommand{\F}{{\mathbb F}}
\renewcommand{\wp}{{\mathfrak p}}
\renewcommand{\P}{{\mathbb P}}
\renewcommand{\O}{{\mathcal O}}
\newcommand{\Pic}{{\rm Pic\,}}
\newcommand{\Ext}{{\rm Ext}\,}
\newcommand{\rank}{{\rm rk}\,}
\newcommand{\sbull}{{\scriptstyle{\bullet}}}
\newcommand{\bX}{X_{\overline{k}}}
\newcommand{\ch}{\operatorname{CH}}
\newcommand{\tors}{\text{tors}}
\newcommand{\cris}{\text{cris}}
\newcommand{\alg}{\text{alg}}
\newcommand{\tX}{{\tilde{X}}}
\newcommand{\tL}{{\tilde{L}}}
\newcommand{\Hom}{{\rm Hom}}
\newcommand{\spec}{{\rm Spec}}
\newcommand{\M}{{\mathfrak{M}}}
\newcommand{\Gm}{{\mathbb{G}_m}}
\let\isom=\simeq
\let\rk=\rank
\let\tensor=\otimes
\newcommand{\X}{\mathfrak{X}}
\newcommand{\mydot}{{\small{\bullet}}}
\newcommand{\NS}{{\rm NS}}
\let\into=\hookrightarrow
\newcommand{\piloc}{{\pi^{loc}}}

\newtheorem{theorem}[equation]{Theorem}      
\newtheorem{lemma}[equation]{Lemma}          %
\newtheorem{corollary}[equation]{Corollary}  
\newtheorem{proposition}[equation]{Proposition}
\newtheorem{scholium}[equation]{Scholium}

\theoremstyle{definition}
\newtheorem{conj}[equation]{Conjecture}
\newtheorem*{example}{Example}
\newtheorem{question}[equation]{Question}

\theoremstyle{definition}
\newtheorem{remark}[equation]{Remark}

\numberwithin{equation}{subsection}

\title{On the Picard group of moduli spaces}
\author{Kirti Joshi}
\address{Math. department, University of Arizona, 617 N Santa Rita, Tucson
85721-0089, USA.} \email{kirti@math.arizona.edu}
\author{V.~B.~Mehta}
\address{School of Mathematics, Tata Institute of Fundamental Research, Homi Bhabha Road,
Mumbai 400 005, India.} \email{vikram@math.tifr.res.in}
\date{Version \jobname.tex of July 11, 2009, processed \today}

\begin{abstract}
We study the Picard groups of moduli spaces in positive
characteristics and we give a ``$p$-adic'' proof that the Picard
group of moduli of vector bundles of fixed determinant is isomorphic
to the group of integers. Along the way we prove that the local
fundamental group scheme of a normal unirational projective variety
is trivial. This is reminiscent of results of Serre and Nygaard who
studied the fundamental groups of smooth, projective, unirational
varieties.
\end{abstract}

\maketitle


\section{Introduction}
\subsection{Notations}
Throughout this paper, we work over an algebraically closed field
$k$ and we will assume that $k$ has characteristic $p>0$ unless
stated otherwise (the letter $p$ will be reserved for the
characteristic of $k$). The letter $\ell$ will denote a prime
distinct from characteristic of $k$ (so $\ell\neq p$).

Let $C/k$ be a smooth, projective curve over $k$ of genus $g\geq 2$.
Let $L\in \Pic(C)$ be a line bundle on $C$, and let $\M(r,L)$ be the
moduli of semistable vector bundles on $C$ of rank $r$ and
determinant $L$.

It is well-known, that if the degree $\deg(L)=d$ is coprime to the
rank $r$, then $\M(r,L)$ is a smooth, projective variety; if the
degree and the rank are not coprime, then $\M(r,d)$ is singular in
general (the exception is genus $g=2,L=\O_C, r=2$).

When $k$ has characteristic zero, the Picard groups of $\M(r,L)$
have been studied in \cite{drezet89} where it has been shown that
$\Pic(\M(r,L))=\Z$. In characteristic $p>0$, one expects this result
but as far as we are aware, the result has not been established in
literature. The proof of \cite{drezet89} does not seem to adapt
readily to positive characteristics because of the failure of
``Kempf's Lemma'' which is crucial in the proof of \cite{drezet89}.

In this note we prove that $\Pic(\M(r,L))=\Z$ in the following cases
(1) the coprime case and (2) in the case when $L=\O_C$. Our proof
is, in some sense, ``a $p$-adic'' proof and has the merit of working
uniformly in all cases and may also work in the case of moduli of
$G$-bundles ($G$ semi-simple) though we have not verified this at
the moment.

As was pointed out to us by Norbert Hoffmann, that the fact that the
Picard group of the moduli space of vector bundles on a curve is
$\Z$ can also be established by using the relationship between
moduli spaces and moduli stacks especially the corresponding results
for moduli stacks (see \cite{beauville98} and its references for
characteristic zero case and \cite{faltings03,biswas08} for positive
characteristic) and the methods of \cite{beauville98}. This is, of
course, of independent interest but we do not pursue this point of
view here. Instead we show how the result on the Picard group may be
established by a study of the local fundamental group scheme of the
moduli space. In the course of this we establish that the local
fundamental group scheme of a normal unirational variety $X$ is
trivial (a result of independent interest) and in particular, if its
\'etale fundamental group is  also trivial, then we deduce that
$H^1(X,\O_X)$ is zero for such a variety. This should be thought of
as a complement to results of \cite{serre59}, \cite{nygaard78} and
\cite{nori82} on the \'etale fundamental group (scheme) of smooth,
unirational varieties.

To summarize our strategy in brief: we prove that the local
fundamental group scheme (which classifies $F$-trivial vector
bundles) is trivial. Using this we prove that $\Pic(\M(r,L))$ is
reduced and discrete. Next we show that $\Pic(\M(r,L))$ is
torsion-free so we reduce to computing the rank of the
N\'eron-Severi group. This in turn is easily reduced to showing that
the second betti number of $\M(r,L)$ is one (in the cases of
interest).

We would like to thank Norbert Hoffmann and Georg Hein for
conversations. This work was done while the first author was
visiting the Tata Institute of Fundamental Research and he  would
like to thank the Institute for its hospitality.

\section{Unirationality of $\M(r,L)$}
\subsection{} Let $X/k$ be a projective variety, we will say that
$X$ is unirational if there is a dominant, \emph{separable} morphism
$\P^n\to X$ for some $n\geq 1$.

\subsection{} In this section we recall that
$\M(r,L)$ is a unirational variety. In characteristic zero this is
proved in \cite{seshadri82} and in positive characteristic this is
proved in \cite{hein08}.

\begin{theorem}
Let $C$ be a smooth, projective curve and let $L$ be a line bundle
on $C$, let $\M(r,L)$ be the moduli space of semistable vector
bundles on $C$ of rank $r$ and with determinant $L$. Then $\M(r,L)$
is a unirational variety.
\end{theorem}

\section{The fundamental group scheme of a normal unirational variety}
\subsection{}      In \cite{serre59}  it was shown that any smooth,
projective unirational variety over an algebraically closed field of
characteristic zero has trivial fundamental group. This result was
subsequently extended to smooth, unirational threefolds in positive
characteristic in \cite{nygaard78}.

   The main result of this section is to prove that if $X$ is a normal unirational variety, then
its local fundamental group scheme is trivial. We recall that the
local fundamental group scheme, denoted here by $\pi^{\rm loc}(X)$,
was introduced in \cite{mehta02}. We also note that when $X$ is a
smooth, unirational variety the result follows from \cite{nori82} as
$\pi^{\rm loc}$ is a quotient of $\pi^{\rm Nori}(X)$ and it was
shown in \cite{nori82} that the latter is trivial.

\begin{theorem}\label{pilocX}
Let $X/k$ be a normal, unirational variety. Then $\pi^{\rm
loc}(X)=1$.
\end{theorem}

\subsection{}
Let us record the following important corollary of
Theorem~\ref{pilocX}.

\begin{corollary}\label{h1vanishing}
Let $X/k$ be a normal, unirational variety. Assume
that $X$ is simply connected. Then $H^1(X,\O_X)=0$.
\end{corollary}
\begin{proof}
As $X$ is simply connected, we see that every \'etale $\Z/p$-cover
is trivial. By \cite{SGA7} the semi-simple part of $F:H^1(\O_X)\to
H^1(\O_X)$ is controlled by \'etale $\Z/p$-covers of $X$ so  the
semi-simple part is trivial by our hypothesis that $X$ is simply
connected. Hence we conclude that the action of Frobenius on
$H^1(X,\O_X)$ is purely nilpotent. Suppose, if possible that,
$H^1(X,\O_X)$ is not trivial. Let $V$ be a non-trivial extension of
$\O_X$ by $\O_X$. As $F$ is nilpotent on $H^1(X,\O_X)$, we see that
$V$ is $F$-trivial. By Theorem~\ref{pilocX} we see that $V$ is
trivial. This is a contradiction as we had assumed that $V$ is
non-trivial. This completes the proof.
\end{proof}

\subsection{}
The proof of Theorem~\ref{pilocX} is somewhat delicate and we  will
be break it up into several steps. Before we begin with the proof,
we need some notational detail. Recall that $X$ is unirational, so
there exists a dominant, separable morphism $\P^n\to X$ for some
$n\geq 1$. Blowing up and normalizing if required, we may assume
that we have a morphisms $Z\to \P^n$ and $Z\to X$ which makes the
diagram commute and the arrow $Z\to Y\to X$ is the Stein
factorization of $Z\to X$.
\begin{equation}\label{maindiag}
\xymatrix{
  Z \ar[d]_{} \ar[r]^{}\ar[dr]
                & Y \ar[d]^{}  \\
  \P^n \ar@{.>}[r]|-{}
                & X             }
\end{equation}

\begin{proposition}\label{piZ}
In the notation of the diagram \ref{maindiag} above we have
$\piloc(Z)=1$.
\end{proposition}
\begin{proof}
Suppose that $V$ is an $F$-trivial vector bundle on $Z$. Then we may
assume that $F^*(V)=\O_Z^r$ for $r=\rank(V)$. Let $B^1_Z$ be the
sheaf of locally exact differentials. Alternately, we may define
$B^1_Z$ be the short exact sequence of sheaves
\begin{equation}\label{eq:b1def}
0\to\O_Z\to F_*(\O_Z)\to B^1_Z\to 0.
\end{equation}
 We claim that  $V\tensor B_Z^1$
has no sections and in fact $V$ has exactly $r=\rank(V)$ sections.

Suppose, if possible, that $H^0(B_Z^1\tensor V)\neq 0$. Then we have
a non-zero map $V^*\to B^1_Z$. But the Harder-Narasimhan flag of
$B^1_Z$ has negative slopes (this is so because the restriction of
$B^1_Z$ to a suitable open set is isomorphic to the restriction of
$B^1_{\P^n}$ to an open set whose complement has codimension at
least two and the latter is easily seen to have negative slopes),
while that of $V^*$ has non-negative. So such a map must be zero.
Thus we have proved that $B^1_Z\tensor V$ has no sections. Now we
show that this implies that $V$ has $r=\rank(V)$ sections.

Tensoring \eqref{eq:b1def}  with $V$ gives
$$0\to V\to F_*(\O_Z)\tensor V\to B^1_Z\tensor V\to 0.$$
Hence we have from the cohomology long exact
sequence:
$$0\to H^0(V)\to H^0(F^*(V))\to H^0(B^1\tensor V)=0$$
and $H^0(F^*(V))=H^0(\O_Z)^r$.

So $V$ has numerically trivial determinant and $H^0(V)=r=\rank(V)$
and $F^*(V)$ is trivial and  such a $V$ is semistable (with respect
to any polarization  of $Z$). By \cite[Lemma 3.4, page~59]{ein80},
such a $V$ is trivial.

This proves the assertion.
\end{proof}

\begin{proposition}\label{piY}
We claim that for $Y$ as in the diagram~\ref{maindiag}, we have
$$\piloc(Y)=1.$$
\end{proposition}
\begin{proof}
Observe that if $f:Z\to Y$ has connected fibres so
$f_*(\O_Z)=\O_Y$. Suppose that $V$ is an $F$-trivial
vector bundle on $Y$. Then $f^*(V)$ is $F$-trivial on
$Z$ so by Proposition~\ref{piZ} we see that
$f^*(V)=\O_Z^r$ for some $r\geq $ and as
$f_*(\O_Z)=\O_Y$ we see that $V=\O_Y^r$. This proves
the assertion.
\end{proof}

\begin{theorem}
Let $X$ be a normal unirational variety over an algebraically closed
field of characteristic $p>0$. Then we have $\piloc(X)=1$.
\end{theorem}
\begin{proof}
We may assume that we have a diagram~\ref{maindiag}. Suppose $V$ is
an $F$-trivial vector bundle on $X$. Suppose $E\to X$ is a finite
group scheme cover over which $V$ is trivial. Now the pull-back, to
$Y$, of $V$ arises from $E\times_XY \to Y$. But as the pull-back of
$V$ to $Y$ is trivial by Proposition~\ref{piY}, so we have a section
$Y\to E\times_XY \to E$, so by composition have a map $Y\to E$. On
the other hand, as $Y$ is reduced, this map factors as $Y\to
E^{red}\to X$. But $Y\to X$ is separable, so $E^{red}\isom X$, but
this gives a section $X\to E^{red} \to E$ which is a contradiction.
\end{proof}

\section{$\M(r,L)$ is simply connected}
\subsection{}
In this section we outline a proof that $\M(r,L)$ is simply
connected. When $k$ has characteristic zero, this is immediate from
\cite{serre59} by the unirationality of $\M(r,L)$.
\begin{proposition}\label{pitrivial}
We have $\pi_{1}^{et}(\M(r,d))=1$, that is, $\M(r,L)$
is simply connected.
\end{proposition}
\begin{proof}
This is proved using the fact that $\M(r,L)$ lives in a family over
$W=W(k)$ the ring of Witt vectors of $k$. Indeed by \cite{balaji08},
we know that there is a moduli scheme $\M(r,L)\to\spec(W)$ whose
special fibre is our $\M(r,L)$. Moreover, in characteristic zero, it
is well-known, that $\M(r,L)$ is simply connected. Now we use the
specialization theorem of \cite[Sect. 2, Corollary 2.4, Expos\'e
X]{sga1} to deduce that $\M(r,L)/k$ is also simply-connected.
\end{proof}

\section{The coprime case}
\subsection{} Let us first consider the case when the degree
$d=\deg(L)$ and the rank $r$ are coprime. In this
case $\M(r,L)$ is a smooth, projective variety, and
hence we will refer to this case (of coprime degree
and rank) as the ``smooth case''. This case also
explains our  strategy of proof. Let us recall what
we want to prove.
\begin{theorem}
Assume that $L$ is a line bundle on $C$ of degree $d$
and that $(r,d)=1$, then we have $\Pic(\M(r,L))=\Z$.
\end{theorem}

\begin{lemma}
Under the hypothesis at the beginning of this
section, we have $\Pic(\M(r,L))$ has no
$\ell$-torsion for $\ell\neq p$.
\end{lemma}
\begin{proof}
This follows from Lemma~\ref{pitrivial}: as $\pi_1^{et}(\M(r,L))=1$,
we see that $\M(r,L)$ does not admit any $\Z/\ell^m$-covers for any
$m\geq 1$. On the other hand, if $\Pic(\M(r,L))$ has $\ell$-torsion
for some $\ell$, then certainly $\M(r,L)$ has \'etale
$\Z/\ell^m$-covers for this $\ell$ and some $m\geq 1$.
\end{proof}

\begin{lemma}
Under the hypothesis at the beginning of this
section, we have the following
\begin{enumerate}
\item $H^1_{dR}(\M(r,L)/k)=0,$
\item $H^2_{cris}(\M(r,L)/W)$ is torsion free,
\item $\Pic(\M(r,L))$ has no $p$-torsion.
\end{enumerate}
\end{lemma}
\begin{proof}
By \cite{illusie79} we have the implications
$(1)\Rightarrow (2)\Rightarrow (3)$. But for the
convenience of the reader we will prove all of them.

Since $\M(r,L)$ is unirational by \cite{serre59}, we
have $H^0(\M,\Omega^1_{\M})=0$ and we have already
established that $H^1(\M,\O_\M)=0$, so we have the
vanishing of the $H^1_{dR}$ as required for (1).

To prove (2) we note that by the universal
coefficient theorem (see \cite{illusie79}) we have
$$0\to H^1_{cris}(\M(r,L)/W)\tensor_W k\to
H^1_{dR}(\M(r,L)/k)\to
Tor_1(H^2_{cris}(\M(r,L)/W),k)\to 0.$$ By (1) the
term in the middle is zero so we are done, and in
particular we have also proved that
$$H^1_{cris}(\M(r,L)/W)=0.$$

From Corollary~\ref{h1vanishing} we know that
$H^1(\M(r,L),\O_{\M(r,L)})=0,$ it follows that
$\Pic(\M(r,L))$ is reduced and discrete. Hence
$\Pic(\M(r,L))\isom \NS(\M(r,L))$. By
\cite{illusie79}, we have
$$ NS(\M(r,L))\tensor W\into H^2_{cris}(\M(r,L)/W),$$
so the former is torsion free follows from the fact
that the latter is torsion free.

Thus we have proved all the assertions.
\end{proof}

\begin{proposition}
Under the hypothesis that the degree and the rank are coprime, we
have that the rank of N\'eron-Severi group of $\M(r,L)$ is one.
\end{proposition}
\begin{proof}
Since $NS(\M(r,d))\into H^2_{cris}(\M(r,L)/W)$, its suffices to
prove that the latter has rank one. We may compute the rank after
tensoring with $\Q_p$. So it suffices to compute the second Betti
number $b_2(\M(r,L))$. By \cite{katz74} this may be computed using
$\ell$-adic cohomology instead of crystalline cohomology.

Now we may further reduce to the case when $k$ is a
finite field. So we will assume for the rest of the
proof that $k$ is a finite field. In this case the
required Betti number may be extracted from the
computation of the Poincar\'e polynomial of $\M(r,L)$
(see \cite{harder74}). We carry out this computation
now.

Following Harder and Narasimhan (see \cite{harder74}), let
$$Z_m(T)=\frac{(1+T^{1-2m})^{2g}}{(1-T^{-2m})(1-T^{2(1-m)})},$$
and let
$$P(T)=T^{2(r^2-1)(g-1)}Z_2(T)\cdots Z_r(T).$$
Then expand $P(T)$ as a power series in $T^{-1}$, say
$$P(T)=T^{2(r^2-1)(g-1)}\sum_{v=0}^\infty b_v T^{-v}.$$
 Then $$b_v=\dim H^v_{et}(\M(r,L),\Q_\ell)$$ for
 $0\leq v\leq 2((r-1)(g-1)+d$ (here $d=\deg(L)$).
So to prove the assertion that $b_2=1$, it suffices
to prove that coefficient of $T^{-2}$ in the
expansion of $Z_2(T)\cdots Z_r(T)$ is one. Now
writing this out we have
$$Z_2(T)\cdots
Z_r(T)=\prod_{j=2}^r\frac{(1+T^{1-2j})^{2g}}{(1-T^{-2j})(1-T^{2(1-j)})}.$$
As $j\geq 2$, the only term which contributes as
$T^{-2}$ in $\sum_v b_v T^{-v}$ is the term
$\frac{1}{(1-T^{2(1-j)})}$ for $j=2$, and as this is
$$\frac{1}{1-T^{-2}}=1+T^{-2}+T^{-4}+\cdots,$$
so this term contributes exactly one to $b_2$ and so
we see that $b_2(\M(r,L))=1$.

This completes our proof in the coprime case
\end{proof}
\section{The trivial determinant case}
\subsection{}
We now consider the case when $L=\O_C$. In this case the
moduli space is singular and we do not know a simple formula for the
Poincar\'e polynomial which is valid in arbitrary characteristics.
Hence to compute the Betti number $b_2(\M(r,\O_C))$ we have to
resort to calculating the Betti number in characteristic zero. This
calculation is well-known, but we indicate a proof for completeness.
\begin{theorem}
Let $C$ be a smooth, projective curve of genus $g\geq 2$ over an
algebraically closed field $k$ of  characteristic $p>0$. Then we
have $\Pic(\M(r,\O_C))=\Z$.
\end{theorem}
\begin{proof}
For simplicity of notation let us write $X=\M(r,\O_C)$ for the
moduli space of semistable bundles of rank $r$ and trivial
determinant. Then we have shown that $H^1(X,\O_X)=0$ and hence we
deduce that $\Pic(X)$ is a discrete scheme. Next we observe that
$\Pic(X)$ has no $p$-torsion or $\ell$-torsion. Let us first
consider $p$-torsion. Then as $\pi^{loc}(X)=1$, we see that $X$ has
no $\alpha_p$ or $\mu_p$ covers. Further as $$H^1(X,\O_X)=0$$ and as
$H^1(X,\Z/p)$ is the semisimple part of the Frobenius map
$H^1(X,\O_X)\to H^1(X,\O_X)$ so $H^1(X,\Z/p)=0$. Hence $X$ has no
\'etale $p$-covers. So $\Pic(X)$ has no $p$-torsion. Next as
$\pi^{et}_1(X)=1$ we see that $X$ does not have $\ell$-covers for
any $\ell\neq p$. Thus $\Pic(X)$ has no $\ell$-torsion. So $\Pic(X)$
is torsion free as claimed. Thus $\Pic(X)=NS(X)$. Next we remark
that by \cite[Page 145]{grothendieck-dix-exposes} that
$NS(X)\tensor\Q_\ell\to H^2_{'et}(X,\Q_\ell(1))$ is injective. This
is a standard argument but we give a proof for completeness. Let
$\Gm$ be the multiplicative group scheme. Then we have an exact
sequence (for each $n\geq 1$) of sheaves for the \'etale topology:
\begin{equation}
1\to \mu_{\ell^n}\to \Gm\to \Gm\to 1,
\end{equation}
This induces a long exact sequence
\begin{equation}
\cdots H^1(X,\mu_{\ell^n})\to \Pic(X)\to\Pic(X)\to
H^2(X,\mu_{\ell^n})\to\cdots
\end{equation}
From this we extract
\begin{equation}
0\to \Pic(X)/\ell^n\to H^2(X,\mu_{\ell^n})\to Br'(X)[\ell^n]\to 0.
\end{equation}
where $Br'(X)=H^2(X,\Gm)$ is the cohomological Brauer group of $X$
and $Br'(X)[\ell^n]$ is its $\ell^n$ torsion. In the inverse limit
we get (see \cite[Page 145]{grothendieck-dix-exposes})
\begin{equation}
0\to NS(X)\tensor\Z_\ell\to H^2(X,\Z_\ell(1))\to T_\ell Br'(X)\to 0,
\end{equation}
where $T_\ell Br'(X)$ is the $\ell$-adic Tate module of $Br'(X)$. At
any rate this shows that the rank of $NS(X)$ is at most  $b_2(X)$
(the second Betti number of $X$); since $X$ is projective the rank
is at least one. So it remains to prove that $b_2(X)=1$. This is the
content of the next proposition.
\end{proof}


\begin{proposition}
Let $C$ be a smooth, projective curve over an algebraically closed
field $k$ of characteristic $p>0$. Let $\M(r,\O_C)$ be the moduli of
semistable vector bundles on $C$ of rank $r$ and trivial
determinant. Then we have
$$b_2(\M(r,\O_C))=\dim
H^2_{et}(\M(r,\O_C),\Q_\ell)=1.$$
\end{proposition}
\begin{proof}
Let $W=W(k)$ be the ring of Witt vectors of $k$, let $K$ be the
quotient field of $W$. Then we have a family $\M\to \spec(W)$ whose
special is $\M(r,\O_C)$ and whose generic fibre is
$\M(r,\O_{C_K})_K$ for some lift $C_K$ of $C$ to $K$. By the proper
base change theorem, to prove our assertion it suffices to prove
that $H^2_{et}(\M(r,\O_{C_K}),\Q_\ell)$ has dimension one. Thus we
reduce to calculating the Betti number in characteristic zero. So we
may assume that $k=\C$ and that $C/k$ is a smooth, projective curve
of genus $g\geq 2$. In this case it is well-known, see for instance
\cite{boysal07}, that we have an isomorphism
$$\Pic(\M(r,\O_C))\isom H^2(\M(r,\O_C),\Z),$$
and hence by the main result of \cite{drezet89} we deduce that the
Betti number $$b_2(\M(r,\O_C))=1.$$ This completes the proof of our
assertion.
\end{proof}


\begin{thebibliography}{{N}or82}

\bibitem[BH08]{biswas08}
I.~Biswas and Norbert Hoffmann, \emph{The line bundles on moduli stacks of
  principal bundles on a curve},
  Preprint:\href{http://arxiv.org/abs/0805.2915v2}{http://arxiv.org/abs/0805.2%
915v2}, 2008.

\bibitem[BK05]{boysal07}
Arzu Boysal and Shrawan Kumar, \emph{Explicit determination of the {P}icard
  group of moduli spaces of semistable {$G$}-bundles on curves}, Math. Ann.
  \textbf{332} (2005), no.~4, 823--842.

\bibitem[BLS98]{beauville98}
Arnaud Beauville, Yves Laszlo, and Christoph Sorger, \emph{The {P}icard group
  of the moduli of {$G$}-bundles on a curve}, Compositio Math. \textbf{112}
  (1998), no.~2, 183--216.

\bibitem[DN89]{drezet89}
J.-M. Drezet and M.~S. Narasimhan, \emph{Groupe de {P}icard des vari\'et\'es de
  modules de fibr\'es semi-stables sur les courbes alg\'ebriques}, Invent.
  Math. \textbf{97} (1989), no.~1, 53--94.

\bibitem[Ein80]{ein80}
Lawrence Ein, \emph{Stable vector bundles on projective spaces in char $p>0$},
  Math. {Ann.} \textbf{254} (1980), 53--72.

\bibitem[Fal03]{faltings03}
Gerd Faltings, \emph{Algebraic loop groups and moduli spaces of bundles}, J.
  Eur. Math. Soc. (JEMS) \textbf{5} (2003), no.~1, 41--68.

\bibitem[Gro68]{grothendieck-dix-exposes}
Alexander Grothendieck, \emph{Le groupe de {B}rauer. {III}. {E}xemples et
  compl\'ements}, Dix {E}xpos\'es sur la {C}ohomologie des {S}ch\'emas,
  North-Holland, Amsterdam, 1968, pp.~88--188.

\bibitem[Hei]{hein08}
Georg Hein, \emph{{$SU(r,L)$} is separably unirational},
  \href{http://arxiv.org/abs/0810.4079}{{http://arxiv.org/abs/0810.4079}}.

\bibitem[HN75]{harder74}
G.~Harder and M.~S. Narasimhan, \emph{On the cohomology groups of moduli spaces
  of vector bundles on curves}, Math. Ann. \textbf{212} (1974/75), 215--248.

\bibitem[Ill79]{illusie79}
Luc Illusie, \emph{Complexe de de\thinspace {R}ham-{W}itt et cohomologie
  cristalline}, Ann. Sci. \'Ecole Norm. Sup. (4) \textbf{12} (1979), no.~4,
  501--661.

\bibitem[KM74]{katz74}
Nicholas~M. Katz and William Messing, \emph{Some consequences of the {R}iemann
  hypothesis for varieties over finite fields}, Invent. Math. \textbf{23}
  (1974), 73--77.

\bibitem[MS02]{mehta02}
V.~B. Mehta and S.~Subramanian, \emph{On the fundamental group scheme}, Invent.
  Math. \textbf{148} (2002), no.~1, 143--150.

\bibitem[{N}or82]{nori82}
Madhav~{V}. {N}ori, \emph{The fundamental group scheme}, Proceedings of
  {I}ndian {A}cademy of {S}cience ({M}ath. {S}ci.) \textbf{91} (1982), no.~2,
  73--122.

\bibitem[Nyg78]{nygaard78}
Niels Nygaard, \emph{On the fundamental group of a unirational {$3$}-fold},
  Invent. Math. \textbf{44} (1978), no.~1, 75--86.

\bibitem[Ser59]{serre59}
J.-P. Serre, \emph{On the fundamental group of a unirational variety}, J.
  London Math. Soc. \textbf{34} (1959), 481--484.

\bibitem[Ses82]{seshadri82}
C.~S. Seshadri, \emph{Fibr\'es vectoriels sur les courbes alg\'ebriques},
  Ast\'erisque, vol.~96, Soci\'et\'e Math\'ematique de France, Paris, 1982,
  Notes written by J.-M. Drezet from a course at the \'Ecole Normale
  Sup\'erieure, June 1980.

\bibitem[sga71]{sga1}
\emph{Rev\^etements \'etales et groupe fondamental}, Springer-Verlag, Berlin,
  1971, S\'eminaire de G\'eom\'etrie Alg\'ebrique du Bois Marie 1960--1961 (SGA
  1), Dirig\'e par Alexandre Grothendieck. Augment\'e de deux expos\'es de M.
  Raynaud, Lecture Notes in Mathematics, Vol. 224.

\bibitem[SGA73]{SGA7}
\emph{Groupes de monodromie en g\'eom\'etrie alg\'ebrique. {II}},
  Springer-Verlag, Berlin, 1973, S\'eminaire de G\'eom\'etrie Alg\'ebrique du
  Bois-Marie 1967--1969 (SGA 7 II), Dirig\'e par P. Deligne et N. Katz, Lecture
  Notes in Mathematics, Vol. 340.

\bibitem[VBM08]{balaji08}
T.~E. {V}enkata {B}alaji and {V}.~{B}. {Mehta}, \emph{Singularities of moduli
  spaces of vector bundles over curves in characteristic $0$ and $p$}, Michigan
  {M}ath. {J}. \textbf{57} (2008), 37--42.

\end{thebibliography}
\providecommand{\bysame}{\leavevmode\hbox to3em{\hrulefill}\thinspace}
\providecommand{\MR}{\relax\ifhmode\unskip\space\fi MR }
\providecommand{\MRhref}[2]{%
  \href{http://www.ams.org/mathscinet-getitem?mr=#1}{#2}
}
\providecommand{\href}[2]{#2}

\end{document}